\documentclass[11pt]{amsart}
\usepackage{amssymb, enumerate, xspace, graphicx} 
\usepackage[all]{xy}
\SelectTips{cm}{} 

\numberwithin{equation}{section}

1

\numberwithin{subsection}{section}

\allowdisplaybreaks[1]

\newtheorem*{namedtheorem}{\theoremname}
\newcommand{\theoremname}{testing}

\newtheorem{theorem}{Theorem}[section] \newtheorem*{theorem*}{Theorem}

\newtheorem{proposition-definition}[theorem] {Proposition-Definition}

\newtheorem{lemma}[theorem]{Lemma}

\theoremstyle{definition}

\newtheorem{remark}[theorem]{Remark} 

\theoremstyle{remark} 

\newcommand\cO{\mathcal{O}}

\renewcommand\AA{\mathbb{A}} 
\newcommand\CC{\mathbb{C}} 
 \newcommand\FF{\mathbb{F}}
\newcommand\GG{\mathbb{G}} 
\newcommand\II{\mathbb{I}}

 \newcommand\PP{\mathbb{P}}
\newcommand\QQ{\mathbb{Q}}

 \newcommand\ZZ{\mathbb{Z}}

\newcommand\rO{\mathrm{O}}

\newcommand\rmc{\mathrm{c}}

\newcommand\rmm{\mathrm{m}}

\newcommand\arr{\ifinner\to\else\longrightarrow\fi}

\newcommand\xarr{\xrightarrow}

\newcommand{\larr}{\longrightarrow}

\renewcommand\H{\operatorname{H}}

\newcommand\eqdef{\overset{\mathrm{\scriptscriptstyle def}} =}

\newcommand\into{\hookrightarrow}

\renewcommand\th{^\text{th}}

\def\displaytimes_#1{\mathrel{\mathop{\times}\limits_{#1}}}

\def\displayotimes_#1{\mathrel{\mathop{\bigotimes}\limits_{#1}}}

\newcommand\spec{\operatorname{Spec}}

\newcommand\generate[1]{\langle #1 \rangle}

\newcommand{\PGL}{\operatorname{PGL}}

\newdir{ >}{{}*!/-5pt/@{>}}

\newcommand\doublelong[2]{\mathbin{\xymatrix{{}\ar@<3pt>[r]^{#1}
\ar@<-3pt>[r]_{#2}&}}}

\newcommand\gm{\GG_{\rmm}} 
\newcommand\A[1][*]{\operatorname{A}^{#1}}

\renewcommand\H[1][*]{\operatorname{H}^{#1}}
\renewcommand\O[1][n]{\rO_{#1}} \newcommand\SO[1][n]{\mathrm{SO}_{#1}}
\newcommand\GL[1][n]{\mathrm{GL}_{#1}}
\newcommand\SL[1][n]{\mathrm{SL}_{#1}} 
\newcommand\Sp[1][n]{\mathrm{Sp}_{#1}} \newcommand{\An}{\AA^n}
\newcommand{\0}{\{0\} } 
\renewcommand\c{\rmc} \newcommand\mmu[1][2]{\boldsymbol{\mu}_{#1}}

\begin{document}

\title[On the Chow rings of classifying spaces for classical groups]{On
the Chow rings of classifying spaces\\for classical groups}

\author[Alberto Molina]{Luis Alberto Molina Rojas}

\address[Alberto Molina]{Dipartimento di Matematica\\ Universit\`a di
Roma Tre\\ Largo San Leonardo Murialdo, 1\\ I-00146 Roma\\ Italy}
\email{molina@mat.uniroma3.it}

\author[Angelo Vistoli]{Angelo Vistoli}

\address[Angelo Vistoli]{Dipartimento di Matematica\\ Universit\`a di
Bologna\\ Piazza di Porta San Donato 5\\ I-40126 Bologna\\ Italy}
\email{vistoli@dm.unibo.it}


\begin{abstract} We show how the stratification method, introduced by Vezzosi in his study of $\PGL_{3}$, provides a unified approach to the known computations
of the Chow rings of the classifying spaces of $\GL$, $\SL$, $\Sp$, $\O$
and $\SO$. \end{abstract}

\maketitle



\section{Introduction}

To algebraic topologists, the cohomology of classifying spaces of linear algebraic groups (or, equivalently, of compact Lie groups) has been an important object of study for a long time. Recently, B{.}~Totaro (\cite{totaro-classifying}) has introduced an algebraic analogue of this cohomology, the Chow ring of the classifying space of a linear algebraic group $G$, denoted by $\A_{G}$. If $\H_{G}$ denotes the integral cohomology of the classifying space of $G$, there is a natural ring homomorphism $\A_{G} \arr \H_{G}$, which is, in general, neither surjective not injective.

Remarkably, the Chow ring $\A_{G}$ seems to be better behaved than $\H_{G}$: for example, if $G$ is a finite abelian group $\A_{G}$ is the symmetric algebra on the group of characters, while $\H_{G}$ is much larger (unless $G$ is cyclic). This is truly surprising to someone who is familiar with the theory of Chow rings of smooth projective algebraic varieties, because these tends to be infinitely more complicated and less computable than their cohomology.

Rationally, the situation is very well understood. If $G$ is a connected
algebraic group, then the homomorphism $\A_{G}\otimes\QQ \arr \H_{G}\otimes\QQ$ is an isomorphism, and both rings coincide with the ring of invariants under the Weyl group in the symmetric algebra of the ring of characters of a maximal torus; this is classical, due to Leray and Borel, in the case of cohomology, and to Edidin and Graham (\cite{edidin-graham-characteristic}) for the Chow ring. Furthermore, this ring of invariants is always a polynomial ring, as was shown by Chevalley. With integral coefficients, the situation is much more subtle.

The Chow ring $\A_{G}$ has been computed for the classical groups $\GL$, $\SL$, $\Sp$, $\O$ or $\SO$, but not for the $\PGL_{n}$ series. The results are as follows. Each of the groups above comes with a tautological representation $V$, of dimension $n$ (or $2n$, in the case of $\Sp$). Every  representation $V$ of an algebraic group $G$ has Chern classes $\c_{i}(V) \in \A[i]_{G}$. When $G$ is a classical group, we denote the Chern classes of the tautological representation simply by $c_{i}$.

Burt Totaro (\cite{totaro-classifying}) and R.~Pandharipande
(\cite{rahul-orthogonal}) described $\A_{G}$ when $G = \GL$, $\SL$,
$\Sp$, $\O$ and $\SO$ when $n$ is odd. We will use the following
notation: if $R$ is a ring, $t_{1}$, \dots,~$t_{n}$ are elements of $R$,
$f_{1}$, \dots,~$f_{r}$ are polynomials in $\ZZ[x_{1}, \dots, x_{n}]$,
we write
   \[
   R = \ZZ[t_{1}, \dots, t_{n}]/ \bigl(f_{1}(t_{1}, \dots,
   t_{n}), \dots, f_{r}(t_{1}, \dots, t_{n})\bigr)
   \]
to indicate the the
ring $R$ is generated by $t_{1}$, \dots,~$t_{n}$, and the kernel of the
evaluation map $\ZZ[x_{1}, \dots, x_{n}] \arr R$ sending $x_{i}$ to
$y_{i}$ is generated by $f_{1}$, \dots,~$f_{r}$. When there are no
$f_{i}$ this means that $R$ is a polynomial ring in the $t_{i}$.

First the case of the special groups.

\begin{theorem*}[B. Totaro]\hfill \begin{enumerate}

\item $\A_{\GL} = \ZZ[c_{1}, \dots, c_{n}]$.

\item $\A_{\SL} = \ZZ[c_{2}, \dots, c_{n}]$.

\item $\A_{\Sp} = \ZZ[c_{2}, c_{4}, \dots, c_{2n}]$. \end{enumerate}
\end{theorem*}

The first two cases follow very easily from the well known description
via generators and relations of the Chow ring of a Grassmannian.

In all three cases, the Chow ring is isomorphic to the cohomology ring.

\begin{theorem*}[R. Pandharipande, B. Totaro]\hfill \begin{enumerate}
\item $\A_{\O} = \ZZ[c_{1}, \dots, c_{n}]/(2c_{\mathrm{odd}})$.

\item If $n$ is odd, then $\A_{\SO} = \ZZ[c_{2}, \dots,
c_{n}]/(2c_{\mathrm{odd}})$. \end{enumerate}

\end{theorem*}

The notation $2c_{\mathrm{odd}}$ means ``all the elements $2c_{i}$ for $i$ odd''.

In these cases the Chow ring is not isomorphic to the cohomology ring: the ring  $\H_{\O}$ was computed independently by Brown (see \cite{brown-On}) and Feshbach (\cite{feshbach-On}); the result is considerably more involved.
The fact that these formulae are so simple is another manifestation of
the tamer nature of $\A_{G}$, as opposed to $\H_{G}$. 

When $n$ is odd, then $\O \simeq \SO \times \mmu$, and this allows to obtain the result for $\SO$ from that for $\O$. When $n$ is even this fails, and the situation is more complicated. Even rationally, the Chern classes of the tautological representation do not generate the Chow ring, or the cohomology. It is well known that when $n = 2m$, the tautological representation has an Euler class $\epsilon_{m} \in \H[2m]_{\SO}$, whose square is $(-1)^{m}c_{n}$: this class, together with the even Chern classes $c_{2}$, $c_{4}$, \dots,~$c_{n-2}$ generate $\A_{\SO}\otimes\QQ = \H_{\SO}\otimes\QQ$. Totaro noticed that when $n = 4$ the class $\epsilon_{2}$ is not in the image of $\A_{\SO}$; shortly afterwards, Edidin and Graham (\cite{edidin-graham-quadric}) constructed a class $y_{m} \in \A[m]_{\SO}$, whose image in $\H_{\SO}$ is, rationally, $2^{m-1}\epsilon_{m}$.

Subsequently, Pandharipande computed $\A_{\SO[4]}$: he showed that it is generated by $c_{2}$, $c_{3}$, $c_{4}$ and $y_{2}$, and gave the relations (his description of the class $y_{2}$ is different, but equivalent to that of Edidin and Graham). Finally, in her Ph.D. thesis Rebecca Field obtained the general result (\cite{field-so2n}), which is as follows.

\begin{theorem*}[R. Field] When $n = 2m$ is even, then
   \[
   \A_{\SO} = \ZZ[c_{2},
   \dots, c_{n}, y_{m}]/ \bigl(y_{m}^{2} - (-1)^{m}2^{n-2}c_{n},\
   2c_{\mathrm{odd}}, \ y_{m}c_{\mathrm{odd}}\bigr).
   \]
\end{theorem*}

Furthermore there are many results, due to Totaro himself (\cite{totaro-classifying}), to P.~Guillot (\cite{guillot-chevalley} and \cite{guillot-steenrod}), and to N.~Yagita (\cite{yagita-extraspecial}), for finite groups. Chow rings of classifiying spaces of exceptional groups have been studied by Yagita (see \cite{yagita-F4}, \cite{yagita-AHss}).

As we already mentioned, the $\PGL_{n}$ series is much harder (this is
an example of a universal phenomenon, that of all the classical groups,
these are the ones giving rise to the deepest problems). For $n = 2$ we have that $\PGL_{2} =
\SO[3]$, and for this group everything is well understood. The cohomology with $\ZZ/3\ZZ$ of the classifying space of $\PGL_{3}$ has been computed coefficients in \cite{kono-mimura-shimada}, and that of $\PGL_{n}$ with $\ZZ/2\ZZ$ coefficients when $n \equiv 2 \pmod{4}$ in \cite{kono-mimura} and \cite{toda-classifying}; furthermore, several results on the cohomology of $\PGL_{p}$ with $\ZZ/p\ZZ$ coefficients were proved in \cite{vavpetic-viruel}. To our knowledge, not much else was known about the cohomology of the classifying space of $\PGL_{n}$ with $\ZZ/p\ZZ$ coefficients, when $p$ divides $n$.

Concerning the Chow ring, for $n = 3$ there is a difficult paper of G{.}~Vezzosi (\cite{vezzosi-pgl3}), where he describes $\A_{\PGL_{3}}$ almost completely. Here is his basic idea. The fundamental tool is the equivariant intersection theory that Edidin and Graham (\cite{edidin-graham-equivariant}) have forged starting from Totaro's idea. Vezzosi stratifies the adjoint representation $\mathfrak{sl}_{3}$ of $\PGL_{3}$ by type of Jordan canonical form, compute the Chow ring of each stratum, and then get generators for $\A_{\PGL_{3}}$ using the localization sequence for equivariant Chow groups. To get relations he restricts to appropriate subgroups of $\PGL_{3}$. His technique has been refined and improved by the second author in \cite{vistoli-pglp}, where he studies the Chow ring and the cohomology of the classifying space of $\PGL_{p}$, where $p$ is an odd prime.

The purpose of this article is to show how this stratification method provides a unified approach to all the known results on the Chow ring of classical groups. Consider a classical group $G$ with its tautological representation $V$. Then one stratifies $V$ in strata in which the stabilizers  are, up to an extension by a unipotent group, smaller classical group. Using the localization sequence for equivariant Chow groups this gives generators for the Chow rings, with relations that come out naturally. To show that the relations suffice, one restricts to appropriate subgroups of $G$ (a maximal torus first, to show that the relations suffice up to torsion, then to some finite subgroup to handle torsion). This turns out to be reasonably straightforward for all the classical groups, except for $\PGL_{n}$. From our point of view, this is due to the fact that the natural representation to use for $\PGL_{n}$, which is the adjoint representation, has a much more complicated orbit structure then in the case of the other groups.

For the cases of $\Sp$ and $\O$, Totaro's proofs, based on his very interesting and important \cite[Proposition~14.2]{totaro-classifying}, are much simpler. In the case of $\SO$ for even $n$, Totaro's method, as implemented by Field, does not seem easier than the stratification method. In the case of $\PGL_{n}$, the stratification method provides the best know results; but there is also a very recent preprint of Kameko and Yagita where they also compute the additive structure of $\A_{\PGL_{n}}$ and $\H_{\PGL_{n}}$, with completely different methods, using the Adams spectral sequence for Brown--Peterson cohomology (\cite{kameko-yagita-pup}).

A few words about the future\footnote{The authors are well aware of the risks involved in making predictions, as people always play Chesterton's game ``Cheat the Prophet''.}. Despite its elementary nature, the stratification methods is powerful; also, as was pointed out to the second author by N{.}~Yagita, it might yield interesting results when applied to more general cohomology theories. However, it seems clear that to proceed much further one will eventually need to introduce substantial amounts of homological machinery, as provided by the theory of motivic cohomology of Voevodsky and Morel. Thus, the way is indicated by the work of N{.}~Yagita (see for example \cite{yagita-modp} and \cite{yagita-AHss}).

\subsection*{Acknowledgments}

Pierre Guillot pointed out an error in the first version of this paper. We thank him heartily.

The second authors would like to acknowledge the very interesting discussions he has had with Nobuaki Yagita and Andrzej Weber on the subject of this paper.

\section{Preliminaries and notation}

In this section we recall some definitions and notations, and state some tecnical results that will be used throughout this paper.

All schemes and algebraic spaces are assumed to be  of finite type over an fixed field $k$. Let $G$ a $g$-dimensional linear algebraic group over $k$, and $X$ a smooth scheme over $k$ with a $G$-action. 

Edidin and Graham (\cite{edidin-graham-equivariant}), expanding on the idea of Totaro, have defined the $G$-equivariant Chow ring of $X$, denoted $\A_{G}(X)$, as follows. For each $i \geq 0$, choose a representation  $V$ of $G$ with an open susbscheme $U\subset V$ on which $G$ acts freely (in which case we call $(V,U)$ a \emph{good pair for $G$}), and such that the codimension of $V\setminus U$ is greater than $i$. The action of $G$ on $X \times U$ is also free, and the quotient $(X\times U)/G$ exists as a smooth  algebraic space; then Edidin and Graham define
   \[
   \A[i]_{G}(X)\eqdef \A[i]\bigl((X\times U)/G\bigr),
   \]
where the right hand term is the usual Chow group of classes of cycles of codimension $i$. This is easily seen to be independent of the good pair $(V,U)$ chosen. Then one sets
   \[
   \A_{G}(X) \eqdef \bigoplus_{i \geq 0}\A[i]_{G}(X).
   \]
If $G$ acts freely on $X$, then there is a quotient $X/G$ as an algebraic space of finite type over $k$, and the projection $X \arr X/G$ makes $X$ into a $G$-torsor over $X/G$; and in this case the ring $\A_{G}(X)$ is canonically isomorphic to $\A(X/G)$

Totaro's definition of the Chow ring of a classifying space is a particular case of this, as
   \[
   \A_{G} = \A_{G}(\spec k).
   \]

The formal properties of ordinary Chow rings extend to equivariant Chow rings. We recall briefly the properties that we need, which will be used without comments in the paper, referring to \cite{edidin-graham-equivariant} for the details.

If $f \colon X \arr Y$ is an equivariant morphism of smooth $G$-schemes there is an induced ring homomorphism $f^{*}\colon \A_{G}(Y) \arr \A_{G}(X)$, making $\A_{G}$ into a contravariant functor from smooth $G$-schemes to graded commutative rings. Furthermore, if $f$ is proper there is an induced homomorphism of groups $f_{*}\colon \A_{G}(X) \arr \A_{G}(Y)$; the projection formula holds.

There is also a functoriality in the group: if $H \arr G$ is a homomorphism of algebraic groups, the action of $G$ on $X$ induces an action of $H$ on $X$, and there is homomorphism of graded rings $\A_{G}(X) \arr \A_{H}(X)$. When $H$ is a subroup of $G$ we will refer to this as a \emph{restriction homomorphism}.

If $H$ is a subgroup of $G$, then there is an $H$-equivariant embedding $X$ into $X\times G/H$, defined in set-theoretic terms by sending $x$ into $(x,1)$. Then the composite of the restriction homomorphism $\A_{G}(X\times G/H) \arr \A_{H}(X\times G/H)$ with the pullback $\A_{H}(X\times G/H) \arr \A_{H}(X)$ is an isomorphism.

Of paramount importance is the localization sequence; if $Y$ is a closed $G$-invariant subscheme of $X$, and we denote by $i\colon Y \into X$ and $j \colon X \setminus Y \into X$ the inclusions, then the sequence
   \[
   \A_{G}(Y) \overset{i_{*}}\larr \A_{G}(X)
   \overset{j^{*}}\larr \A_{G}(X \setminus Y) \arr 0
   \]
is exact.

Furthermore, if $E$ is a $G$-equivariant vector bundle on $X$, there are Chern classes $\c_{i}(E) \in \A[i]_{G}(X)$, enjoying the usual properties. Also, the pullback $\A_{G}(X) \arr \A_{G}(E)$ is an isomorphism.

In particular, since the equivariant vector bundles over $\spec k$ are the representations of $G$, we get Chern classes $\c(V) \in \A[i]_{G}$ for every representation of $G$; and the pullback $\A_{G} \arr \A_{G}(V)$ is an isomorphism.

We also need other easy properties of equivariant Chow rings, for which we do not have a suitable reference.

\begin{lemma}\label{lem:subgroup-free} Let $G$ a linear algebraic group, $X$ a smooth $G$-scheme, $H$ a normal algebraic subgroup $G$. Suppose that the action of $H$ on $X$ is free with quotient $X/H$. Then there is canonical isomorphism of graded rings
   \[
   \A_G(X) \simeq \A_{G/H}(X/H).
   \]
\end{lemma} 

\begin{proof} Let $(V,U)$ be a good pair for $G$, such that the codimension of $V \setminus U$ is greater then $i$. Then 
\begin{align*}
   \A[i]_{G}(X) &= \A[i]\bigl((X\times U)/G\bigr)\\
   &=\A[i] \bigl( ((X\times U)/H)/(G/H)\bigr)\\
   &=\A[i]_{G/H}\bigl((X\times U)/H\bigr).
\end{align*}

Now, the quotient $(X \times V)/H$ is a $G/H$-equivariant vector bundle over $X/H$, $(X\times U)/H$ is an open subscheme of $(X \times V)/H$ whose complement has codimension larger than $i$. This yields isomorphisms
\begin{align*}
   \A[i]_{G/H}\bigl((X\times U)/H\bigr) &\simeq 
   \A[i]_{G/H}\bigl((X\times V)/H\bigr)\\
   &\simeq \A[i]_{G/H}(X/H).
\end{align*}
The resulting isomorphisms $\A[i]_{G}(X) \simeq \A[i]_{G/H}(X/H)$ yield the desired ring isomorphism $\A_{G}(X) \simeq \A_{G/H}(X/H)$.
\end{proof}

\begin{lemma}\label{lem:describe-comp-0}
Let $G$ be an affine linear
group acting on a smooth scheme $X$, $E \arr X$ an equivariant vector bundle of rank $r$. Call $E_{0} \subseteq E$ the complement of the zero section of $E$. Then the pullback homomorphism $\A_{G}(X) \arr \A_{G}(E_{0})$ is surjective, and its kernel is generated by the top Chern class $\c_{r}(E) \in \A[r]_{G}(X)$.
\end{lemma}

\begin{proof}
Call $s\colon X \arr E$ the zero-section. Then the statement follows immediately from the exactness of the localization sequence
   \[
   \A_G(X) \overset{s_{*}}{\longrightarrow} \A_G(E) \arr \A_G(E_0)\arr 0,
   \]
from the fact that the pullback $s^{*}\colon \A_{G}(E) \arr \A_{G}(X)$ is an isomorphism, and from the self-intersection formula, which implies that the composite $\A_G(X) \xarr{s_{*}} \A_G(E) \xarr{s^{*}} \A_G(X)$ is multiplication by $\c_r(E)$.
\end{proof}

\begin{lemma}\label{lem:unipotent-ext}
Let $H$ a linear algebraic group with an isomorphism $H\simeq \AA^n_{k}$ of varieties, such that the for any field extension $k \subseteq k'$ and any $h \in H(k')$, the action of $h$ on $H_{k'}$ by multiplication is corresponds to an affine automorphism of $\AA^{n}_{k'}$ under the isomorphism above. Furthermore, let $G$ be a linear algebraic group acting on $H$ via group automorphisms, that corresponds to a linear action of $G$ on $\AA^{n}_{k}$ under this isomorphism. 

If $G$ acts on a smooth scheme $X$: form the semidirect product $G \ltimes H$, and let $G \ltimes H$ act on $X$ via the projection $G \ltimes H \arr G$. Then the homomorphism
   \[
   \A_{G}(X) \arr \A_{G \ltimes H}(X)
   \]
induced by the projection $G \ltimes H \arr G$ is an isomorphism.
\end{lemma}

\begin{proof} Let $(V,U)$ (resp. $(V',U')$) be a good pair for $G\ltimes H$ (resp. $G$). Then $G\ltimes H$ acts on $U'$ via the projection $G\ltimes H\rightarrow G$: it follows that $G\ltimes H$ acts on $X\times H\times U\times U'$, and since the action of $G\ltimes H$ on $H$ is transitive, and the stabilizer of the origin is $H$, there is an isomorphism
\begin{align*}
   (X\times H \times U\times U')/(G\ltimes H) &= (X\times (G\ltimes H)/H \times U\times U')/(G\ltimes H)\\
   &\simeq (X\times U\times U')/G.
\end{align*}
Look at the following commutative diagram:
   \[
   \xymatrix{ (X\times H \times U \times U')/(G\ltimes H)  \ar[r] \ar[d]^{\pi_1} 
   &(X\times U \times U')/G \ar[d]^{\pi_2}  \\
   (X\times U\times U')/(G\ltimes H) \ar[r]^-{f} & (X\times U')/G .}
   \]
Note that $\pi_1$ is an affine bundle: in fact, it is a fiber bundle with fiber isomorphic to $\AA^n$, and structure group $G\ltimes H$ that acts on $\AA^n$ by affine transformations, since the action of $G$ on $H$ is affine and the action of $H$ on itself is affine. It follows from \cite[p{.}~35]{grothendieck-chevalley} that $\pi_1^*$ is an isomorphism. On the other hand, since $U\times U'$ is an open set of $V\times V'$ on which $G$ acts freely, $\pi_2^*$ is the identity on the equivariant Chow ring $\A_G(X)$, up to a degree that can be made arbitrarily large: so we have a commutative triangle 
   \[
   \xymatrix{  & \A_{G\ltimes H}(X\times H) \ar@{<-}[ld]_{\pi_1^*}^{\simeq}
   \ar@{<-}[rd]_{\simeq} &   \\ \A_{G\ltimes H}(X)  & & \A_G(X) \ar[ll]_{f^*} }
   \]
where the horizontal arrow is exactly the map induced by the projection $G\ltimes H\rightarrow G$.
\end{proof}

Here is another auxiliary result: it is well known (see for instance \cite{vezzosi-pgl3}) that $\A_{\mmu[n]}\simeq \ZZ [\xi ]/(n\xi )$, where $\xi$ is the first Chern class of the character given by the inclusion $\mmu[n] \hookrightarrow \gm$. If $G$ is an algebraic group, we will denote by $\xi \in \A_{G\times \mmu[n]}$ the image of $\xi$ under the map $\A_{\mmu[n]} \rightarrow \A_{G\times \mmu[n]}$ induced by the projection $G\times \mmu[n]\rightarrow \mmu[n]$. Using the projection $G\times \mmu[n]\rightarrow G$, we can consider $\A_{G\times \mmu[n]}$ as an $\A_G$-algebra. Then $\A_{G\times \mmu[n]}$ admits the following description:

\begin{lemma}\label{lem:product-cyclic}
As an $\A_G$-algebra, $\A_{G\times \mmu[n]}$ is generated by the element $\xi$, and the kernel of the evaluation map $\A_G[x]\rightarrow \A_G[\xi]$ is the ideal $(n x)$. In other words,
   \[
   \A_{G\times \mmu[n]} = \A_G[\xi ]/(n\xi ).
   \] 
\end{lemma} 

\begin{proof} The action of $\mmu[n]$ on $\AA^1$ given by the embedding $\mmu[n] \hookrightarrow \gm$ can be extended to an action of $G \times \mmu[n]$ by letting $G$ act trivially on $\AA^1$. Then from Lemma~\ref{lem:describe-comp-0} we have that $\A_{G\times \mmu[n]} \rightarrow \A_{G\times \mmu[n]}(\gm )$ is surjective, and its kernel is generated by $\xi$. Since $\gm / \mmu[n] \simeq \gm$, from Lemma~\ref{lem:subgroup-free} we deduce that $\A_{G\times \mmu[n]}(\gm )\simeq \A_G(\gm )$, and since $G$ acts trivially on $\gm$ and $\gm$ is an open subset of the affine line, $\A_G(\gm )\simeq \A_G$. So we have that $\A_{G\times \mmu[n]} \simeq \A_G [\xi ]/(n\xi )$, as claimed. \end{proof}

\section{The special groups: $\GL$, $\SL$ and $\Sp$}

Let us fix a field $k$: we write $\GL$, $\SL$ and $\Sp$ for the corresponding algebraic groups over $k$.

These groups are always much easier to study: they are special, in the sense that every principal bundle is Zariski locally trivial. For $\GL$ and $\Sp$ the idea works in a very similar way: let us work out $\Sp$, that is marginally harder. We proceed by induction on $n$, the case $n=0$ being trivial.

Consider $V = \AA^{2n}$, the tautological representation of $\Sp$, with its symplectic form $h\colon V \times V \arr k$ given in coordinates by
   \[
   h(x_{1}, \dots, x_{2n}, y_{1}, \dots, y_{2n}) = x_{1}y_{n+1} + \dots +
   x_{n}y_{2n} - x_{n+1}y_{1} - \dots -x_{2n}y_{n}.
   \]
Denote by $e_{1}$, \dots,~$e_{2n}$ the canonical basis of $V$.

The orbit structure of $V$ is very simple: there are two orbits, the origin and its complement $U \eqdef V \setminus\{0\}$. Consider the subspace
   \[
   V' = \generate{e_{1}, \dots, e_{n-1}, e_{n+1}, \dots,
   e_{2n-1}};
   \]
the restriction of $h$ to $V'$ is a non-degenerate symplectic form, and $V  = V' \oplus \generate{e_{n}, e_{2n}}$. This induces an embedding $\Sp[n-1] \hookrightarrow \Sp$, identifying $\Sp[n-1]$ with the stabilizer of the pair $(e_n,e_{2n})$. 

Let $G$ the stabilizer of the element $e_n$: then we have that $\Sp[n-1]\subseteq G\subseteq \Sp[n]$. The first inclusion admits a splitting: if $A\in G$, then $A$ stabilizes the orthogonal complement $\left\langle e_n\right\rangle ^{\perp}$. It follows that $A$ induces a linear endomorphism on the quotient $\left\langle e_n\right\rangle ^{\perp}/\left\langle e_n\right\rangle \simeq V'$, and this endomorphism is easily seen to preserve the symplectic form $h|_{V'}$, so it is an element of $\Sp[n-1]$. Thus we have a projection $G\rightarrow \Sp[n-1]$: let $H$ its kernel, so that $G=\Sp[n-1] \ltimes H$.

The structure of $H$ is as follows; the matrices in $H$ are exactly those for which there are scalars $a_{1}$, \dots,~$a_{2n-1}$ such that
   \[
   Ae_{i} =
   \begin{cases}
   e_{i} - a_{i+n}e_{n} & \text{if $i = 1$, \dots,~$n-1$}\\
   e_{n} &\text{if $i = n$}\\
   e_{i} + a_{i-n}e_{n} & \text{if $i = n+1$, \dots,~$2n-1$}\\
   a_{1}e_{2} + \dots + a_{2n-1}e_{2n-1} + e_{2n} &\text{if $i = 2n$}.\\
   \end{cases}
   \]
This yields an isomorphism of varieties $H\simeq \AA ^{2n-1}$. It is not hard to see that the conditions of Lemma~\ref{lem:unipotent-ext} are satisfied for the action of  $\Sp[n-1]$ on $H$; hence the embedding $\Sp[n-1] \subseteq G$ induces an isomorphism of rings $\A_{G} \simeq
\A_{\Sp[n-1]}$, so the composite
   \[
   \A_{\Sp}(U) \arr \A_{\Sp[n-1]}(U) \arr \A_{\Sp[n-1]}(e_{n}) =
   \A_{\Sp[n-1]}
   \]
is an isomorphism. The restriction of the representation $V$ to $\Sp[n-1]$ is the direct sum of $V'$ and of a trivial 2-dimensional representation: hence the Chern classes $c_{i} = \c_{i}(V)$ restrict to the $\c_{i}(V')$. From the induction hypothesis, we conclude that $\A_{\Sp}(U)$ is generated by the images of $c_{2}$, \dots,~$c_{2n-2}$.

From Lemma~\ref{lem:describe-comp-0} we conclude that every class in $\A_{\Sp}$ can be written as a polynomial in $c_{2}$, \dots,~$c_{2n-2}$, plus a multiple of $c_{2n}$. By induction on the degree we conclude that $c_{2}$, \dots,~$c_{2n}$ generate $\A_{\Sp}$.

To prove their algebraic independence, let us restrict to $\A_{T_{n}}$, where $T_{n} \simeq \gm^{n}$ is the standard maximal torus in $\Sp$, consisting of diagonal matrices with entries $(t_{1}, \dots, t_{n}, t_{1}^{-1}, \dots, t_{n}^{-1})$. Then $\A_{T_{n}}$ is the polynomial ring $\ZZ[\xi_{1}, \dots, \xi_{n}]$, where $\xi_{i}$ is the first Chern class of the 1-dimensional representation given by the $i\th$ projection $T_{n} \arr \gm$. Then the total Chern class of the restriction of $V_{n}$ to $T_{n}$ is
   \[
   (1+\xi_{1}) \dots (1+\xi_{n}) (1-\xi_{1}) \dots
   (1-\xi_{n}) = (1-\xi_{1}^{2}) \dots (1-\xi_{n}^{2});
   \]
hence the restriction of $c_{2i}$ is the $i\th$ elementary symmetric function of $-\xi_{1}^{2}$, \dots,~$-\xi_{n}^{2}$. This proves the independence of the $c_{2i}$.

As we mentioned, the argument for $\GL$ is very similar. For $\SL$, one can proceed similarly, but it is easier to use the fact that, if $\GL$ acts freely on an algebraic variety $U$, the induced morphism $U/\SL \arr U/\GL$ makes $U/\SL$ into a principal $\gm$-bundle on $U/\GL$, associated with the determinant homomorphism $\det\colon \GL \arr \gm$. Hence, by Lemma~\ref{lem:describe-comp-0}, we have an isomorphism $\A_{\SL} \simeq \A_{\GL}/(c_1)$, which gives us what we want.

\begin{remark}
All these arguments work with cohomology. when $k = \CC$. The localization sequence in cohomology does not quite work in the same way, as the restriction homomorphism from the cohomology of the total space to that of an open subset is not necessarily surjective. However, if $Y$ is a smooth closed subvariety of a smooth complex algebraic variety $X$, of pure codimension $d$, then there is an exact sequence
   \[
   \cdots \arr \H[i-2d]_{G}(Y) \arr \H[i]_{G}(X) \arr
   \H[i]_{G}(X\setminus Y) \arr \H[i-2d+1]_{G}(Y) \arr \cdots.
   \]
Hence if we know that either the pullback $\H_{G}(X) \arr \H_{G}(X\setminus Y)$ is surjective, or the pushforward $\H_{G}(Y) \arr \H_{G}(X)$ is injective, we can conclude that we have an exact sequence
   \[
   0 \arr
   \H_{G}(Y) \arr \H_{G}(X) \arr \H_{G}(X\setminus Y) \arr 0;
   \]
and this is sufficient to mimic the arguments above and give the result for cohomology.
\end{remark}

\begin{remark}
These results can also be proved very simply from a result of Edidin and Graham (see \cite{edidin-graham-characteristic}): if $G$ is a special algebraic group, $T$ a maximal torus and $W$ the Weyl group, the natural restriction homomorphism $\A_{G} \arr \bigl(\A_{T}\bigr)^{W}$ is an isomorphism.
\end{remark}

\section{The Chow ring of the classifying space of $\O$}

Let us fix a field $k$ of characteristic different from $2$. If $V = k^{n}$ is an $n$-dimensional vector space, we define a quadratic form $q\colon V \arr k$ in the standard form
   \[
   q(x_{1}, \dots, x_{n}) =
   x_1x_{m+1}+\dots +x_mx_{2m}
   \]
when $n = 2m$, and
   \[
   q(x_{1}, \dots,
   x_{n})=x_1x_{m+1}+\dots +x_mx_{2m}+ x_{2m+1}^{2}
   \]
when $n = 2m+1$. We will denote by $\O$ the algebraic group of linear transformations preserving this quadratic form.

\begin{theorem}[R. Pandharipande, B. Totaro]
\[
\A_{\O} =
\ZZ[c_1,\dots,c_n]/(2c_\mathrm{odd}).
\]
\end{theorem}

\begin{remark}\label{rmk:bitorsor} Let $V'$ be another $n$-dimensional vector space over $k$, with a non-degenerate quadratic form $q'\colon V' \arr k$. We can associate with this another algebraic group $\rO(q')$, which will not be isomorphic to $\O = \rO(q)$, in general, unless $k$ is algebraically closed.

However, one can show that there is an isomorphism of Chow rings $\A_{\O} \simeq \A_{\rO(q')}$, such that the classes $\c_{i}(V)$ in the left hand side correspond to the classes $\c_{i}(V')$ in the right hand side. The principle that allows to prove this has been known for a long time (\cite{giraud}): it is the existence of a bitorsor $I \arr \spec k$. This is the scheme representing  the functor of isomorphisms of $(V,q)$ with $(V', q')$. On $I$ there is a left action of $\rO(q')$ and right action of $\O$, by composition. These two actions commute, and make $I$ into a torsor under both groups (because $(V,q)$ and $(V',q')$ become isomorphic after a base extensions).

In general, assume that $G$ and $G'$ are algebraic groups over a field $k$ (in fact, any algebraic space will do as a base), and $I \arr \spec k$ is $(G', G)$-bitorsor: that is, on $I$ there is a right action of $G$ and left action of $G'$, and this makes $I$ into a torsor under both groups. If $X$ is a $k$-algebraic space on which $G'$ acts on the left, then we can produce a $k$-algebraic space $I\times^{G} X$ on which $G$ acts on the left, by dividing the product $I\times_{\spec k} X$ by the right action of $G$, defined by the usual formula $(i, x)g = (ig, xg^{-1})$. The left action of $G'$ is by multiplication on the first component: the quotients $G \backslash X$ and $G' \backslash (I\times^{G} X)$ are canonically isomorphic.

This operation gives an equivalence of the category of $G$-algebraic spaces with the category of $G'$-algebraic spaces. When applied to representations, it yields representations, and gives an equivalence of the category of representations of $G$ and of $G'$. Furthermore, given a representation $V$ of $G$, with an open subset $U \subseteq V$ on which $G$ acts freely, we get a representation $V' = I\times^{G} V$ with an open subset $U' = I\times^{G} U$ on which $G'$ acts freely, so that the quotients $G \backslash U$ and $G' \backslash U'$ are isomorphic. In Totaro's construction this gives an isomorphism of $\A_{G}$ with $\A_{G'}$.

So, in particular, the result that we have stated for $\O$ also holds for $\rO(q')$ for any other non-degenerate $n$-dimensional quadratic form $q'$, and we have
   \[
   \A_{\rO(q')} =
   \ZZ[c_1,\dots,c_n]/(2c_\mathrm{odd}).
   \]
\end{remark}

The proof of the Theorem will be split into two parts: first we show that  the $c_{i}$ generate $\A_{\O}$, then that ideal of relations is generated by the given ones.

For the first part we proceed by induction on $n$.

For $n=1$, $q(x)=x_1^2$, and $\O[1] = \mmu$, so
   \[
   \A_{\O[1]} =
   \A_{\mmu} \simeq \ZZ[c_1]/(2c_1).
   \]

For $n>1$, let $B=\{v\in \AA^n \mid q(v)\neq 0\}$, and set $Q = q^{-1}(1)$. Then $q\colon B\rightarrow \gm$ is a fibration, with fibers isomorphic to $Q$. This fibration is not trivial, but it becomes trivial after an \'etale base change. Set
   \[
   \widetilde{B} = \{(t,v) \in \gm
   \times B \mid t^2=q(v)\},
   \]
and consider the cartesian diagram
   \[
   \xymatrix{ \widetilde{B} \ar[r] \ar[d] & B \ar[d]^{q} \\
   \gm \ar[r]^{(-)^2} & \gm }
   \]
where the first column is projection onto the first factor, and the top row is defined by the formula $(t,v) \mapsto tv$.

There are obvious commuting actions of $\mmu$ and $\O$ on $\widetilde{B}$, the first defined by $\epsilon\cdot (t,v) = (\epsilon t, v)$, and the second by $M\cdot (t,v) = (t, Mv)$. The quotient $\widetilde{B}/\mmu$ is isomorphic to $B$, and the induced action of $\O$ on the quotient coincides with the given action on $B$. From Lemma~\ref{lem:subgroup-free}, we obtain an isomorphism
   \[
   \A_{\O}(B)\simeq A_{\mmu\times \O}(\widetilde{B}).
   \]

Then there is an isomorphism of $\gm$-schemes $\widetilde{B}\simeq \gm \times Q$ defined by the formula $(t,v) \mapsto (t, v/t)$. The given actions of $\mmu$ and of $\O$ on $\widetilde{B}$ induce commuting actions on $\gm \times Q$ given by $\epsilon \cdot (t,v)=(\epsilon t,\epsilon v)$ for $\epsilon\in \mmu$ and $M(t,v)=(t,Mv)$ for $M\in \O$. These define an action of $\mmu\times \O$ on $\gm \times Q$, and $\A_{\O}(B)$ is isomorphic to $\A_{\mmu\times \O}(\gm \times Q)$.

This action of $\mmu\times \O$ on $\gm \times Q$ extends uniquely to an action of $\mmu\times \O$ on $\AA^{1}\times Q$, defined by the same formulae. This action is defined by two separate action on $\AA^{1}$ and $Q$, and the action on $\AA^{1}$ is linear, defined by the non-trivial character of $\mmu$ through the projection $\mmu \times \O \arr \mmu$. Call $\xi$ the first Chern class of this representation. From Lemma~\ref{lem:describe-comp-0}, we have an isomorphism \begin{equation}\label{formula1} \A_{\mmu\times \O}(\gm \times Q) \simeq \A_{\mmu\times \O}(Q)/(\xi). \end{equation}

To investigate $\A_{\mmu\times \O}(Q)$ we will also use an orthogonal basis
$e'_{1}$, \dots,~$e'_{n}$ of $V$, in which $q$ has the form
   \[
   q(x_{1}e'_{1} + \dots + x_{n}e'_{n}) = x_1^{2} + \dots + x_{m}^{2} -
   x_{m+1}^{2} - \dots - x_{n}^{2}
   \]
when $n = 2m$, and
   \[
   q(x_{1}e'_{1} +
   \dots + x_{n}e'_{n}) = x_1^{2} + \dots + x_{m+1}^{2} - x_{m+2}^{2} -
   \dots - x_{n}^{2}
   \]
when $n = 2m+1$.

Now, the action of $\mmu\times \O$ on $Q$ is transitive; let $H$ the
stabilizer of the point $e'_{1} \in Q$. The structure of $H$ is as
follows. Set $V' \eqdef \generate{e'_{2}, \dots, e'_{n}}$, so that $V$ is
the orthogonal sum $\generate{e'_{1}} \oplus V'$, and call $q'$ the
restriction of $q$ to $V'$. Then the group $\rO_{q'}$ of linear
automorphisms of $V'$ preserving $q'$ is naturally embedded into $\O$,
as the stabilizer of $e'_{1}$. Notice that in an appropriate basis $q'$
has the standard form
   \[
   q'(x_{1}, \dots, x_{n-1}) = x_1x_{m+1}+\dots
   +x_mx_{2m}
   \]
when $n = 2m+1$, and the opposite of the standard form
   \[
   q'(x_{1}, \dots, x_{n-1}) = -(x_1x_{m}+\dots +x_{m-1}x_{2m-2}+
   x_{2m-1}^{2})
   \]
when $n = 2m$; in both cases the orthogonal group
$\rO(q')$ is isomorphic to $\O[n-1]$, and we identify it with $\O[n-1]$.

The stabilizer of $e'_{1}$ in $\mmu \times \O$ is the group $\mmu \times
\O[n-1]$, embedded into $\mmu \times \O$ with the injective homomorphism
   \[
   (\epsilon, M) \longmapsto (\epsilon, \epsilon M).
   \]
It follows that
   \begin{align*}
   \A_{\mmu \times \O}(Q)
      &\simeq \A_{\mmu \times \O}\bigl((\mmu \times
   \O)/(\mmu \times \O[n-1])\bigr)\\
   &\simeq \A_{\mmu \times \O[n-1]}.
   \end{align*}
We obtain a chain of isomorphisms
   \begin{align*}
   \A_{\O}(B) &\simeq \A_{\mmu\times \O}(Q)/(\xi)\\
   &\A_{\mmu \times \O[n-1]}/(\xi).
   \end{align*}   
Finally, from Lemma~\ref{lem:subgroup-free} we get
an isomorphism
   \begin{align*}
   \A_{\mmu \times \O[n-1]}/(\xi) &\simeq 
   \A_{\O[n-1]}[\xi]/(\xi)\\
   &\simeq \A_{\O[n-1]}.
\end{align*}
The composite $\A_{\O}
\arr \A_{\O}(U) \arr \A_{\O[n-1]}$ is the pullback induced by the
embedding $\O[n-1] \subseteq \O$. 

The restriction of $V$ to $\O[n-1]$ is the direct sum of $V'$ and a
trivial 1-dimensional representation, hence the restriction $\A_{\O}
\arr \A_{\O[n-1]}$ carries $c_{i}$ into $\c_{i}(V')$. Therefore, by
induction hypothesis, the images of $c_{1}$, \dots,~$c_{n-1}$ generate
$\A_{\O}(B)$.

Next, we claim that the restriction homomorphism $\A_{\O}(\AA^{n}
\setminus \{0\}) \arr \A_{\O}(B)$ is an isomorphism. To see this, set
   \[
   C = \{ v\in \AA^{n} \setminus \{0\} \mid q(v)=0\}
   \]
with its reduced
scheme structure, and consider the fundamental exact sequence
   \[
   \A_{\O}(C) \stackrel{i_*}{\longrightarrow} \A_{\O}(\AA^n\setminus \{0\})\arr
   \A_{\O}(U) \arr 0.
   \]
We need to show that $i_*$ is the zero map. In fact, $q\colon \AA^{n} \setminus \{0\} \arr \AA^{1}$ is smooth, since the characteristic of the base field is not $2$, so $C$ is the scheme-theoretic inverse image of $\{0\}$. The map $q\colon \AA^{n} \setminus \{0\} \arr \AA^{1}$ is $\O$-equivariant, if we let $\O$ act
trivially on $\AA^{1}$; and the fundamental class $[0] \in
\A_{\O}(\AA^{1})$ equals zero. Since the inverse image of $[0]$ in
$\A_{\O}(\AA^n \setminus \{0\})$ is $[C]$, we can conclude that
   \[
   [C] = 0 \in
   \A_{\O}(\AA^n \setminus \{0\}).
   \]
Next we show that the pullback $i^{*}\colon
\A_{\O}(\AA^n\setminus \{0\}) \arr \A_{\O}(C)$ is surjective: in this case, for every $\alpha
\in \A_{\O}(C\setminus \{ 0\})$, we have $\alpha =i^*\beta$ for some $\beta \in
\A_{\O}(\AA^n \setminus \{0\})$, so
   \[
   i_*(\alpha)=i_*i^*(\beta )=[C]\cdot \beta =0
   \]
by the projection formula, and $i_*$ is the zero map, as claimed.

To show surjectivity, notice that the action of $\O$ on $C$ is transitive. Let us investigate the stabilizer $G$ of $e_{1}\in C$. Set $n = 2m$ or $n = 2m+1$, as usual. If we define
   \[
   V' = \generate{e_{2},
   \dots, e_{m}, e_{m+2}, \dots, e_{n}}
   \]
then the restriction of $q$ to $V'$ has the standard form, and $V$ is the orthogonal sum $V' \oplus \generate{e_{1}, e_{m+1}}$. This gives an embedding $\O[n-2] \subseteq \O$, identifying $\O[n-2]$ with the stabilizer of the pair $(e_{1}, e_{m+1})$.

An analysis very similar to that we have carried out for the stabilizer of a vector under $\Sp$ leads to the conclusion that the stabilizer $G$ of $e_1$ is a semidirect product $\O[n-2]\ltimes H$, where $H$  is isomorphic to $\AA ^{n-1}$ as a variety, the action of an element of $H$ is itself is given by an affine map, and the action of $\O[n-2]$ on $H$ is linear: by Lemma~\ref{lem:unipotent-ext}, the embedding $\O[n-2] \subseteq G$ induces an isomorphism of rings $\A_{G} \simeq
\A_{\O[n-2]}$, so the composite
   \[
   \A_{\O}(C) \arr
   \A_{\O[n-2]}(C) \arr \A_{\O[n-2]}(e_{1}) = \A_{\O[n-2]}
   \]
is an isomorphism. But the $c_{i}$ restrict in $\A_{\O[n-2]}$ to the Chern
classes of $V'$: hence, by inductions hypothesis, they generate $\A_{\O[n-2]}$. Hence the pullback $\A_{\O} \arr \A_{\O}(C)$ is surjective, as claimed. This ends the proof that the $c_{i}$ generate $\A_{\O}$. Let us investigate the relations.

The quadratic form $q$ induces an isomorphism $V \simeq V^{\vee}$ of representations of $\O$, hence for each $i$ we have $\c_{i}(V) = (-1)^{i}\c_{i}(V)$. This shows that $2c_{i}= 0$ when $i$ is odd.

To show that these generate the ideal of relations among the $c_{i}$, let $J \subseteq \ZZ[x_{1}, \dots, x_{n}]$ be the ideal generated by $2x_{1}$, $2x_{3}$, \dots. Let $P \in \ZZ[x_{1}, \dots, x_{n}]$ be a homogeneous polynomial such that $P(c_{1}, \dots, c_{n}) = 0 \in \A_{\O}$: we need to check that $P$ is in $J$. By modifying $P$ by an element of $J$, we may assume that $P$ is of the form $Q + R$, where $Q$ is a polynomial in the even $x_{i}$, while $R$ is a polynomial in which every monomial contains some $x_{i}$ with $i$ odd, and all of whose coefficients are either $0$ or $1$.

Let $T_{m} \simeq \gm^{m}$ be the standard torus in $\O$: the embedding $T_{m} \subseteq \O$ sends $(t_{1}, \dots, t_{m})$ into the diagonal matrix with entries $(t_{1}, \dots, t_{m}, t_{1}^{-1}, \dots, t_{m}^{-1})$ if $n =2m$, and $(t_{1}, \dots, t_{m}, t_{1}^{-1}, \dots, t_{m}^{-1}, 1)$ if $n = 2m+1$. Then $\A_{T_{m}} = \ZZ[\xi_{1}, \dots, \xi_{m}]$, where $\xi_{i}$ is the first Chern class of the $i\th$ projection $\chi_{i}\colon T_{n} \arr \gm$. The restriction of $V$ to $T_{n}$ splits as $\rho \eqdef \chi_{1} + \dots + \chi_{m} + \chi_{1}^{-1} + \dots + \chi_{m}^{-1}$ when $m$ is even, and $\rho + 1$ when $n$ is odd. Hence the total Chern class of the restriction of $V$ to $T_{n}$ is
   \[
   (1+\xi_{1}) \dots (1+\xi_{m}) (1-\xi_{1}) \dots
   (1-\xi_{m}) = (1-\xi_{1}^{2}) \dots (1-\xi_{m}^{2});
   \]
and this means
that the restrictions of the $c_{i}$ is $0$ when $i$ is odd, while $c_{2j}$ restricts to the $j\th$ symmetric function of $-\xi_{1}^{2}$, \dots,~$-\xi_{m}^{2}$. Hence the restrictions of even Chern classes are algebraically independent. In the decomposition $0 = P(c_{1}, \dots, c_{n}) = Q(c_{2}, \dots, c_{2m}) + R(c_{1}, \dots, c_{n})$ the summand $R(c_{1}, \dots, c_{n})$ restricts to $0$, so $Q(c_{2}, \dots, c_{2m})$ also restricts to $0$. This implies that $Q = 0$. So we have that $P$ has coefficients that are either $0$ or $1$.

Now take a basis $e'_{1}$, \dots,~$e'_{n}$ of $V$ in which $q$ has a diagonal form. Consider the subgroup $\mmu^{n} \subseteq \O$ consisting of linear transformations that take each $e'_{i}$ into $e'_{i}$ or $-e'_{i}$. If we call $\eta_{i}$ the first Chern class of the character obtained composing the $i\th$ projection $\mmu^{n} \arr \mmu$ with the embedding $\mmu \into \gm$, then by Lemma~\ref{lem:product-cyclic} we have
   \[
   \A_{\mmu^{n}} =
   \ZZ[\eta_{1}, \dots, \eta_{n}] /(2\eta_{1}, \dots, 2\eta_{n}).
   \]
There
is a natural ring homomorphism from $\A_{\mmu^{n}}$ into the polymomial ring $\FF_{2}[y_{1}, \dots, y_{n}]$ that sends each $\eta_{i}$ into $y_{i}$. The restriction of $V$ to $\mmu^{n}$ has total Chern class $(1+\eta_{1}) \dots (1+\eta_{n})$; hence the image of $c_{i}$ in $\FF_{2}[y_{1}, \dots, y_{n}]$ is the $i\th$ elementary symmetric polynomial $s_{i}$ in the $y_{i}$. The $s_{i}$ are algebraically independent in $\FF_{2}[y_{1}, \dots, y_{n}]$, the image of $0 = P(c_{1}, \dots , c_{n})$ is $P(s_{1}, \dots, s_{n})$, and $P$ has coefficients that either $0$ or $1$. This implies that $P = 0$, and completes the proof of the theorem.

\section{The Chow ring of the classifying space of $\SO$}

Let $k$ be a field of characteristic different from $2$, set $V = k^{n}$, and let $q\colon V \arr k$ be the same quadratic form as in the previous section. Consider the subgroup $\SO \subseteq \O$ of orthogonal linear transformations of determinant $1$.

If $n$ is odd, $\A_{\SO}$ can be easily computed from $\A_{\O}$, as was noticed in \cite{rahul-orthogonal} and \cite{totaro-classifying}.

\begin{theorem}[R. Pandharipande, B. Totaro] If $n$ is odd, then
   \[
   \A_{\SO} = \ZZ[c_2,\dots,c_n]/(2c_{\mathrm{odd}}=0).
   \]
\end{theorem}

\begin{proof}
When is $n$ odd there is an isomorphism $\O\simeq \mmu\times \SO$; the determinant character $\det\colon \O \arr \mmu$ (whose first Chern class in $\A_{\O}$ is $c_{1}$) corresponds to the projection $\mmu\times \SO \arr \mmu$. Then from Lemma~\ref{lem:product-cyclic} we get that
   \[
   \A_{\SO} \simeq \A_{\O}/(c_{1})
   \]
and the conclusion follows.
\end{proof}

\subsection{The Edidin--Graham construction} From now on we shall assume that $n$ is even, and write $n = 2m$.

In this case, $\A_{\SO}$ is not generated by the Chern classes of the standard representation, not even rationally. This can be seen easily for $n = 2$. We have that $\SO[2]$ consists of matrices of the form
   \[
   \begin{pmatrix} t & 0    \\
   0 & t^{-1}\\ \end{pmatrix}
   \]
and so is isomorphic to $\gm$. Then
   \[
   \A_{\SO[2]} = \A_{\gm} = \ZZ[\xi ],
   \]
where $\xi$ is the first Chern class of the tautological representation $L = \AA^{1}$, on which $\gm$ acts via multiplication. Hence $V = L \oplus L^{\vee}$, so $\c_{2}(V) = -\xi^{2}$.

For general $n$, the vector space $V$ will still split as the direct sum of two totally isotropic subspaces, one dual to the other: however, when $n > 2$ this
splitting is not unique, and the totally isotropic subspaces are not invariant under the action of $\SO$, so $V$ is not a direct sum of two nontrivial representations (and $V$ is in fact irreducible). Still, in topology $V$ has an Euler class $\epsilon_{m} \in \H[2m]_{\SO}$, whose square is $(-1)^{m}c_{m}$. Let us recall Edidin and Graham's construction of an algebraic multiple of $\epsilon_{m}$ (see
\cite{edidin-graham-quadric}).

In what follows we will use the classical conventions for projectivizations and Grassmannians; those seem a little more natural in intersection theory than Grothendieck's. So, if $W$ is a vector space, we denote by $\PP(W)$ the vector space of lines in $W$, and by $\GG(r, W)$ the Grassmannian of subspaces of dimension $r$; and similarly for vector bundles.

Denote by $\II(m,V)$ the smooth subvariety of $\GG(m, V)$ consisting of maximal totally isotropic subspaces of $V$. It is well known that $\O$ acts transitively on $\II(m,V)$, and that $\II(m,V)$ has two connected components, each of which is an orbit under the action of $\SO$. Let us choose one of the orbits, for example, the one containing the subspace $\generate{e_{1}, \dots, e_{m}}$. Every totally isotropic subspace of dimension $m-1$ of $V$ is contained in exactly two maximal totally isotropic subspaces, one in each connected component.

There is a well known equivalence of categories between $\O$-torsors and vector bundles of rank $n$ with a non-degenerate quadratic form. If $E$ is a vector bundle on a scheme $X$ with a non-degenerate quadratic form, this corresponds to a $\O$-torsor $\pi\colon P \arr X$, the torsor of isometries between $E$ and $V \times X$; with this torsor we can associate a $\mmu$-torsor (that is, an \'etale double cover) $P/\SO \arr X$ via the determinant homomorphism $\det\colon \O \arr \mmu$. This cover can be described geometrically as follows.

Consider the subscheme $\II(m,E)$ of totally isotropic subbundles in the relative Grassmannian $\GG(m,E) \arr X$ ; the projection $\II(m, E) \arr X$ is proper and smooth, and each of its geometric fibers has two connected components. Let $\II(m, E) \arr \widetilde{\II}(m,E) \arr X$ be the Stein factorization; then $\widetilde{\II}(m,E) \arr X$ is an \'etale double cover, and is precisely the double cover $P/\SO \arr X$. This can be seen as follows.

On $P$ we have, by definition, an isometry of $\pi^{*}E$ with $V\times P$. In  $V\times P$ we have a maximal totally isotropic subbundle $\generate{e_{1}, \dots, e_{m}} \times P$, so we get a maximal totally isotropic subbundle of $\pi^{*}E$. This defines a morphism $P \arr \II(m,E)$ over $X$; the composite $P \arr  \II(m,E) \arr \widetilde{\II}(m,E)$ induces the desired isomorphism $P/\SO \simeq \widetilde{\II}(m,E)$.

Hence, to give a reduction of structure group of $P \arr X$ to $\SO$ is
equivalent to assigning a section $X \arr \widetilde{\II}(m,E)$. This
gives an equivalence of the groupoid of $\SO$-torsors on $X$ with the
groupoid of vector bundles $E \arr X$ of rank $n$ with a non-degenerated
quadratic form, and a section $X \arr \widetilde{\II}(m,E)$. We shall
refer to such a structure as an \emph{$\SO$-structure on $E$}.

Furthermore, given an $\SO$-structure on $E$, if $f\colon T \arr X$ is a
morphism of algebraic varieties, and $L$ is a totally isotropic
subbundle of $f^{*}E$ of rank $m$, we say that $L$ is \emph{admissible}
if the image of $T$ under the morphism $T \arr \II(m,X)$ corresponding
to $L$ is contained in the inverse image of the given embedding $X
\subseteq \widetilde{\II}(m,E)$.

Here is the construction of Edidin and Graham. We will follow their
notation. Let $E$ be a vector bundle of rank $n$ with an $\SO$-structure
on a smooth algebraic variety $X$. For each $i = 1$, \dots,~$m$ consider
the flag variety $f_{i}\colon Q_{i} \arr X$ of totally isotropic flags
$L_1 \subseteq L_2 \subseteq \dots \subseteq L_{m-i} \subseteq E$, with
each $L_s$ of rank $s$. For each $i$, denote by $L_1 \subseteq L_2
\subseteq \dots \subseteq L_{m-i} \subseteq f_{i}^*E$ the universal flag
on $Q_{i}$. The restriction of the quadratic form to $L_{m-i}^{\perp}$
is degenerate, with radical equal to $L_{m-i}$; hence on $Q_{i}$ there
lives a vector bundle $E_{i} \eqdef L_{m-i}^{\perp}/L_{m-i}$ of rank $2i$
with a non-degenerate quadratic form. For each $i = 1$, \dots,~$m-1$ we
have a projection $\pi_{i}\colon Q_{i-1}\arr Q_{i}$, obtained by
dropping the last totally isotropic subbundle in the chain; and
$Q_{i-1}$ is canonically isomorphic, as a scheme over $Q_{i}$, to the
smooth quadric bundle in $\PP(E_{i})$ defined by the quadratic form on
$E_{i}$. This means that $Q_{i-1}$ is a family of quadrics of dimension
$2(i-1)$ over $Q_{i}$. Let us denote by $h_i \in \A[1](Q_{i-1})$ the
restriction to $Q_{i-1}$ of the class
$\c_1\bigl(\cO_{\PP(E_{i})}(1)\bigr) \in \A[1]\bigl(\PP(E_{i})\bigr)$.

Each bundle $E_{i}$ has a canonical $\SO[n-2i]$-structure. Call $\pi_{i}\colon L_{m-i}^{\perp} \arr E_{i}$ the projection. From each totally isotropic vector subbundle $L \subseteq E_{i}$ of rank $m - i$, we get a totally isotropic vector subbundle $\pi_{i}^{*}L \subseteq L_{m-i}^{\perp} \subseteq f^{*}_{i}E$ of rank $m$; then $L$ is admissible if and only if $\pi_{i}^{*}L$ is admissible.

The universal flag $L_1 \subseteq L_2 \subseteq \dots \subseteq L_{m-1}
\subseteq f_{1}^*E$ on $Q_{1}$ can be completed in a unique way to a
maximal totally isotropic flag $L_{1} \subseteq \dots \subseteq L_{m-1} \subseteq L_{m} \subseteq f_{1}^{*}E$ in such a way that $L_{m}$ is admissible.
Then Edidin and Graham define
   \[
   y_{m}(E) =
   f_{*}\bigl(s\cdot\c_{m}(L_{m})\bigr) \in \A[m](X)
   \]
where we have set
   \[
   s = h_2^2 h_3^4 \dots h_m^{2m-2} \in \A(Q_{1}).
   \]

\begin{remark} In this formula each of the classes $h_{i}$ should be
pulled back to $Q_{1}$. Here, and in what follows, we use the following
convention: when $f\colon Y \arr X$ is a morphism of smooth varieties,
and $\xi \in \A(X)$, we will also write $\xi$ for $f^{*}\xi \in \A(Y)$. Similarly, if $E \rightarrow X$ is a vector bundle, we will also write $E$ for $f^*E$. 
This has the advantage of considerably simplifying notation, and should
not lead to confusion. With this notation, when $f$ is proper the
projection formula reads: if $\xi \in \A(X)$ and $\eta \in \A(Y)$, then
   \[
   f_{*}(\xi\eta) = \xi f_{*}\eta.
   \]
\end{remark}

There is also an inductive definition of $y_{m}(E)$. If $m = 1$ then there
is precisely one totally admissible isotropic line subbundle of $E$, and we have $y_{1}(E) = \c_{1}(L)$, by definition.

For $m > 1$ we have a vector bundle $E_{m-1}$ on $Q_{m-1}$ with an
$\SO[n-2]$-structure.

\begin{lemma} The formula
   \[
   y_{m}(E) =
   -{f_{m-1}}_{*}\bigl(h_{m}^{2m-1}y_{m-1}(E_{m-1})\bigr)
   \]
holds.
\end{lemma}

\begin{proof} To prove this, call $g \colon Q_{1} \arr Q_{m-1}$ the
projection: on $Q_{1}$ we have a flag
   \[
   L_{2}/L_{1} \subseteq
   L_{3}/L_{1} \subseteq \dots \subseteq L_{m-1}/L_{1} \subseteq
   g^{*}E_{m-1}
   \]
that makes $Q_{1}$ into the variety of totally isotropic
flags of length $m-2$ in $E_{m-1}$; we complete this to a maximal
totally isotropic flag by adding $L_{m}/L_{1}$. So we get
   \[
   y_{m-1}(E_{m-1}) = g_{*}\bigl(h_2^2 h_3^4 \dots h_{m-1}^{2m-4}
   \c_{m-1}(L_{m}/L_{1})\bigr).
   \]
On the other hand, on $Q_{m-1} \subseteq
\PP(E)$, the line bundle $L_{1} \subseteq f_{m-1}^{*}E$ is the pullback
of the tautological bundle $\cO_{\PP(E)}(-1)$, so $\c_{1}(L_{1}) =
-h_{m}$. Hence we have
   \[
   \c_{m}(L_{m}) = -h_{m}\c_{m-1}(L_{m}/L_{m-1})
   \]
and
\begin{align*}
   -{f_{m-1}}_{*}\bigl(h_{m}^{2m-1}y_{m-1}(E_{m-1})\bigr) &=
   -{f_{m-1}}_{*}\bigl(h_{m}^{2m-1} g_{*}\bigl(h_2^2 h_3^4 \dots
   h_{m-1}^{2m-4}\\ &\qquad \c_{m-1}(L_{m}/L_{1})\bigr)\bigr) \\ &=
   -{f_{1}}_{*}\bigl(h_2^2 h_3^4 \dots h_{m-1}^{2m-4} h_{m}^{2m-1}
   \c_{m-1}(L_{m}/L_{1})\bigr) \\ &= {f_{1}}_{*}\bigl(h_2^2 h_3^4 \dots
   h_{m-1}^{2m-4} h_{m}^{2m-2} \c_{m}(L_{m})\bigr)\\ &= y_{m}(E)
\end{align*}
as claimed. \end{proof}

The Edidin--Graham class $y_{m} \in \A[m]_{\SO}$ is defined as follows.
Take a representation $W$ of $\SO$ with an open subset $U$ on which
$\SO$ acts freely, and whose complement has codimension larger than $m$.
Call $E$ the vector bundle with an $\SO$-structure associated with the
$\SO$-torsor $U \arr U/\SO$. Then we set
   \[
   y_{m} = y_{m}(E) \in
   \A[m](U/\SO) = \A[m]_{\SO}.
   \]
It is easy to verify that this is
independent of the $W$ and $U$ chosen.

\subsection{The main result}

\begin{theorem*}[R. Field] If $n = 2m$, then
   \[
   \A_{\SO} = \ZZ[c_{2},
   \dots, c_{n}, y_{m}]/ \bigl(y_{m}^{2} - (-1)^{m}2^{n-2}c_{n},\
   2c_{\mathrm{odd}}, \ y_{m}c_{\mathrm{odd}}\bigr).
   \]
\end{theorem*}

\begin{remark} Once again, this result can be extended to other
quadratic forms (compare with Remark~\ref{rmk:bitorsor}). Let $V'$ be
another $n$-dimensional vector space over $k$, with a non-degenerate
quadratic form $q'\colon V' \arr k$. This induces a non-degenerate
quadratic form on the exterior powers $\bigwedge^{i}V'$. Let us assume
that there is an isometry $\bigwedge^{n}V \simeq \bigwedge^{n}V'$.

This is equivalent to the following more concrete condition. We will
write $\det q' \in k^{*}/{k^{*}}^{2}$ for the class in
$k^{*}/{k^{*}}^{2}$ of the determinant of a matrix representing $q'$ in
some basis. Then two $n$-dimensional quadratic forms have isomorphic top
exterior powers if and only if they have the same determinant. Hence the
condition above is equivalent to the equality
   \[ \det q' = (-1)^{m} \in
   k^{*}/{k^{*}}^{2}.
   \]

Fix an isometry $\bigwedge^{n}V \simeq \bigwedge^{n}V'$. We can construct an $\bigl(\mathrm{SO}(q'), \SO\bigr)$-bitorsor $I \arr \spec k$, as the scheme representing the functor of isometries $V \simeq V'$ inducing the fixed isometry $\bigwedge^{n}V \simeq \bigwedge^{n}V'$. So we deduce the following result: if the condition above is satisfied, there exists a class $y_{m} \in \A[m]_{\mathrm{SO}(q')}$, such that
   \[
   \A_{\mathrm{SO}(q')} = \ZZ[c_{2}, \dots, c_{n}, y_{m}]/
   \bigl(y_{m}^{2} - (-1)^{m}2^{n-2}c_{n},\ 2c_{\mathrm{odd}}, \
   y_{m}c_{\mathrm{odd}}\bigr).
   \]
\end{remark}

The proof of the theorem will be split into three parts: first we verify that the classes $c_{i}$ and $y_{m}$ generate $\A_{\SO}$, next that the relations holds, and finally that they generate the ideal of relations.

\subsubsection*{Step 1: The generators} We proceed by induction on $m$. In the case $m = 1$ the statement says that
   \[
   \A_{\SO[1]} = \ZZ[c_{2},
   y_{1}]/(y_{1}^{2} + c_{2}) = \ZZ[y_{1}]
   \]
we have seen that $\SO[1] = \gm$, that $y_{1}$ is the first Chern class of the identity character on $\gm$, and that $c_{2} = -y_{1}^{2}$.

Suppose $m > 1$. Set $B = \{x \in \AA^{n} \mid q(x) \neq 0\}$ and $C = \{x \in \AA^{n} \setminus \{0\} \mid q(x) = 0\}$. Proceeding precisely as for $\O$, one establishes the following results.

\begin{enumerate}

\item Let $e'_{1}$, \dots,~$e'_{n}$ be an orthogonal basis of $V$ in which $q$ has the form
   \[
   q(x_{1}e'_{1} + \dots + x_{n}e'_{n}) = x_1^{2}
   + \dots + x_{m}^{2} - x_{m+1}^{2} - \dots - x_{n}^{2}.
   \]
Then the
stabilizer of $e'_{1}\in B$ in $\SO$ is isomorphic to $\SO[n-1]$, and
the composite
   \[
   \A_{\SO}(B) \arr \A_{\SO[n-1]}(B) \arr
   \A_{\SO[n-1]}(e'_{1}) = \A_{\SO[n-1]}
   \]
is an isomorphism.

\item The stabilizer of the pair $(e_{1}, e_{m+1})$ is isomorphic to $\SO[n-2]$. The composite
   \[
   \A_{\SO}(C) \arr \A_{\SO[n-2]}(C) \arr
   \A_{\SO[n-2]}(e_{1}) = \A_{\SO[n-2]}
   \]
is an isomorphism.
\end{enumerate}

Call $i \colon C \subseteq \AA^{n}  \setminus \{0\}$ and $j \colon B \subseteq \AA^{n}  \setminus \{0\}$ the inclusions. Then
we have an exact sequence
   \[
   \A_{\SO} (C)
   \stackrel{i_{*}}\longrightarrow \A_{\SO} (\An \setminus \0 )
   \stackrel{j^{*}}\longrightarrow
   \A_{\SO}(B) \arr 0.
   \]

By induction hypothesis, we have that $\A_{\SO} (C) \simeq \A_{\SO[n-2]}$ is generated as a ring by $c_{2}$, \dots,~$c_{n-2}$ and $y_{m-1}$. From this, and from the relation $y_{m-1}^{2} - (-1)^{m-1}2^{n-4}c_{n-2}$, we see that $\A_{\SO}(C)$ is generated as a module over $\A_{\SO}$ by $1$ and $y_{m-1}$; hence, since $i_{*}$ is a homomorphism of $\A_{\SO}$-modules, by the projection formula, we see that the kernel of the pullback $\A_{\SO}(\AA^{n} \setminus \0) \arr \A_{\SO}(B)$ is generated as an ideal by $i_{*}1 = [C]$ and $i_{*}y_{m-1}$.

As in the case of $\O$, we see that the fundamental class $[C] \in \A_{\SO}(\AA^{n} \setminus \{0\})$ is $0$, because $C$ is the scheme-theoretic zero-locus of the invariant function $q$. Furthermore, the images of $c_{2}$, \dots~$c_{n-1}$ generate $\A_{\SO}(U) \simeq \A_{\SO[n-1]}$: and this implies that $c_{2}$, \dots,~$c_{n-1}$, together with $i_{*}y_{m-1}$, generate $\A_{\SO}(\AA^{n} \setminus \0) = \A_{\SO}/(c_{n})$. Hence $c_{2}$, \dots,~$c_{n}$, $i_{*}y_{m-1}$ generate $\A_{\SO}$. Next, we have a Lemma.

\begin{lemma}\label{lem:identify-class}
   \[
   i_*y_{m-1} = -y_m \in \A_{\SO}(\An \setminus \0 ).
   \]
\end{lemma} 

\begin{proof} Let $W$ be a representation of $\SO$, and $U$ an open set of $W$ on which the action of $\SO$ is free, and such that the codimension of $W \setminus U$ in $W$ is larger than $m$. The vector bundle associated with the $\SO$-torsor $U \arr U/\SO$ is $E \eqdef (\AA^{n} \times U)/\SO$. We set $X \eqdef \bigl((\AA^{n} \setminus \{0\})\times U\bigr)/{\SO}$, so that $X \subseteq E$ is the complement of the zero section, while $Y \eqdef (C\times  U)/{\SO} \subseteq X$ is the closed subscheme consisting of non-zero isotropic vectors, and $Z\eqdef X \setminus Y$. By a slight abuse of notation, we will denote $i\colon Y\hookrightarrow X$ and $j\colon Z\hookrightarrow X$ the inclusions. Note that there is a tautological section $s\colon X\rightarrow E$ defined set-theoretically by $[u,x]\mapsto [u,x,x]$.

Let us first prove that $j^*y_m=0\in \A_{\SO}(B)$. In fact, the tautological section restricted to $Z$ has the property that $q(s(x))\neq 0$ for all $x$, and so $j^*y_m(E)=y_m(j^*E)=0$, due to the following result.

\begin{lemma}
Let $(E,q) \arr X$ be a rank $n=2m$ vector bundle with a non-degenerate quadratic form. Suppose that there exists a section $s\colon X\longrightarrow E$ such that $q(s(x))\neq 0$ for all $x\in X$. Then $y_m(E)=0$.
\end{lemma}

\begin{proof}
Pulling back to the flag variety $Q_{1} \arr X$, it suffices to show that if $L\subset E$ is a rank $m$ totally isotropic subbundle, then $\c_m(L)=0$. The quadratic form gives a perfect pairing $L \times E/L \longrightarrow \mathcal{O}_X$, so $L^{\vee}\simeq E/L$. On the other hand the line subbundle $\generate{s}$ generated by $s$ has intersection with $L$ equal to $0$ at every point of $X$; hence the composite $\cO_{X} \xarr{w} E \arr E/L$ gives a nowhere vanishing section of $E/L$, so that
   \[
   \c_{m}(L) = (-1)^{m} \c_{m}(E/L) = 0
   \]
as claimed.
\end{proof}

It follows that $y_m=d\cdot i_*y_{m-1}$ with $d\in \ZZ$. We will compute $d$ by restricting to a maximal torus; but first observe that since $\SO[n-2]$ is included in $\SO$ as the stabilizer of the pair $(e_1,e_{m+1})$, there is an isomorphism
\begin{align*}
   (\AA^n \times  U)/{\SO[n-2]} & \longrightarrow \AA^2 \times
      \bigl((\AA^{n-2}\times U)/{\SO[n-2]}\bigr) \\
   [(x_1,\dots ,x_n),u] & \longmapsto \bigl((x_1,x_{m+1}),
   [(x_2,\dots ,x_m,x_{m+2},\dots ,x_n),u]\bigr),
\end{align*}
and that $y_{m-1}\in \A_{\SO[n-2]}$ is the Edidin-Graham class of the vector bundle $\bigl(\AA^{n-2}\times U\bigr)/\SO[n-2] \rightarrow U/\SO[n-2]$.

Now, let $T_{m}\subset \SO$ is, as before, the torus of diagonal matrices with
diagonal entries $t_1$, \dots,~,$t_m$, $t_1^{-1}$, \dots~, $t_m^{-1}$,
and $\xi_{i}$ is the first Chern class of the $i\th$ projection $T_{m}
\arr \gm$.

\begin{lemma}\label{lem:restrict-torus} The formulae
\begin{align*}
   \c_n & = (-1)^m\xi_1^2\dots \xi_m^2 \\
\intertext{and}
   y_m & = 2^{m-1}\xi_1\dots \xi_m
\end{align*}
hold in $\A_{T_{m}} = \ZZ[\xi_{1}, \dots , \xi_{m}]$. 
\end{lemma}

\begin{proof}
Reducing the structure group to $T_m$, the vector bundle $E$ on $U/T_{m}$ associated with the standard representation $T_{m} \into \SO \into \GL$ splits into a direct sum of line bundles $\Lambda_1 \oplus \dots \oplus\Lambda_{2m}$, where the $i\th$ summand is the subbundle associates with the 1-dimensional subspace $\generate{e_{i}} \subseteq V$. For each $i = 1$ , \dots,~$m$ we have $\Lambda_{i+n} \simeq \Lambda_{i}^{\vee}$. Then $E$ has an admissible maximal totally isotropic subbundle $\Lambda_1 \oplus \dots \oplus\Lambda_m$, which pulls back to an admissible totally isotropic subbundle on $Q_{1}$. The first Chern class of $\Lambda_{i}$ in $\A[1](U/T_{m}) = \A[1]_{T_{m}}$ is $\xi_{i}$, for $i = 1$, \dots,~$m$, hence
   \[
   \c_m(\Lambda_1\oplus \dots \oplus \Lambda_m)
   = \xi_{1} \dots \xi_{m}\ \in \A[m]_{T_{m}}
   \]
On the other hand, the top Chern classes of any two admissible totally isotropic subbundles of $Q_{1}$ are the same, by \cite[Theorem 1]{edidin-graham-quadric}, so
\begin{align*}
   y_m &=
   f_*\bigl(s\cdot \c_m(\Lambda_1\oplus \dots \oplus \Lambda_m)\bigr)\\
   &= (f_{*}s)\xi_1\dots \xi_m;
\end{align*}
and it is easy to verify that $f_{*}s = 2^{m-1}$.
\end{proof}

It follows that
   \[
    (\AA^{n-2}\times U)/T_{m-1}
   = \Lambda_2 \oplus \dots \oplus\Lambda_m\oplus \Lambda_2^{\vee} \oplus \dots
   \oplus\Lambda_m^{\vee};
   \]
moreover, since
   \[
   \bigl(U\times(\AA^n\setminus \0)\bigr)/{T_m}
   =(\Lambda_1\oplus \dots \oplus \Lambda_m \oplus \Lambda_1^{\vee}\oplus \dots   
   \oplus \Lambda_m^{\vee})\setminus \0,
   \]
we have 
\begin{align*}
      \A(X) &= \A_{T_m}/(c_n)\\
   &=\ZZ[\xi_1,\dots ,\xi_m]/(\xi_1^2\dots \xi_m^2)
\end{align*}
and our aim is to verify that the equation
\begin{equation}\label{formula2}
i_*y_{m-1}=-2^{m-1}\xi_1 \dots \xi_m
\end{equation}
holds in $\ZZ[\xi_1,\dots ,\xi_m]/(\xi_i^2\dots \xi_m^2)$. 

The inclusion of schemes on $U/T_m$ 
   \[
   (\Lambda_1\oplus \Lambda_1^{\vee})\setminus \0\hookrightarrow
   (\Lambda_1\oplus \dots \oplus \Lambda_m \oplus \Lambda_1^{\vee}\oplus \dots
   \oplus \Lambda_m^{\vee})\setminus \0
   \] 
induces a surjection of rings 
   \[
   \ZZ[\xi_1,\dots ,\xi_m]/(\xi_1^2\dots \xi_m^2)
   \rightarrow  \ZZ[\xi_1,\dots ,\xi_m]/(\xi_1^2);
   \] 
since $\ZZ \xi_1\dots \xi_m$ has trivial intersection with the kernel of this map, we can restrict to $(\Lambda_1\oplus \Lambda_1^{\vee})\setminus \0$ to verify equation \ref{formula2}. There is a cartesian diagram
   \[
   \xymatrix{ (\Lambda_1\setminus \0)
   \sqcup (\Lambda_1^{\vee}\setminus \0) \ar[d] \ar[r]                  &
   (\Lambda_1\oplus \Lambda_1^{\vee})\setminus \0 \ar[d]\\
   Y \ar[r]^i        &    X    } 
   \] 
We set
\begin{align*}
   X' & =(\Lambda_1\oplus \Lambda_1^{\vee})\setminus \0, \\
   \intertext{and}
   Y' & =Y'_1\sqcup Y'_2\\
   &=(\Lambda_1\setminus \0)\sqcup (\Lambda_1^{\vee}\setminus \0);
\end{align*}
call $i'\colon Y' \into X'$ the inclusion.

Also, form the vector bundle on $Y'$ defined as
\begin{align*}
      F &\eqdef \Lambda_2 \oplus \dots \oplus\Lambda_m\oplus \Lambda_2^{\vee} 
   \oplus\dots \oplus\Lambda_m^{\vee}\\
   &= \left\langle s(Y') \right\rangle^{\bot}/\left\langle s(Y') \right\rangle.
\end{align*}
We need to check that
   \[
   i'_{*}y_{m-1}(F) = -2^{m-1}\xi_{1} \dots \xi_{m} \in \A(X').
   \]
For $l = 1$, $2$, call $i'_l\colon Y'_l\hookrightarrow X'$ the inclusion, $s'_l \colon Y'_l\rightarrow i_l'^*E$ the tautological section, $F_l$ the restriction of $F$ to $Y'_{l}$.

Observe that the bundle $\Lambda_{2} \oplus \dots \oplus\Lambda_{m}$ of $F$ is totally isotropic: however, its inverse image in $E$ is $\Lambda_{1} \oplus \dots \Lambda_{m}$ is $\Lambda_{2} \oplus \Lambda_{2} \oplus \dots \oplus\Lambda_{m}$ on $Y_{1}$, but is $\Lambda_{2} \oplus \dots \oplus\Lambda_{m} \oplus \Lambda_{1}^{\vee}$ on $Y_{2}$. The first bundle is admissible, the second one is not. Hence we have
\begin{align*}
   y_{m-1}(F_{1}) = 2^{m-2} \xi_{2} \dots \xi_{m} \in \A(Y'_{1})
\intertext{and}
   y_{m-1}(F_{2}) = -2^{m-2} \xi_{2} \dots \xi_{m} \in \A(Y'_{1}).
\end{align*}
Since we also have $[Y_{1}] = -\xi_{1}$ and $[Y_{2}] = \xi_{1}$ in $\A(X')$, we get
\begin{align*}
   i_*y_{m-1} & = {i_1}_*y_{m-1}(F_{1})+{i_2}_*y_{m-1}(F_{2}) \\
                            &= {i_1}_* i_1^*2^{m-2}\xi_2\dots \xi_m 
                            -{i_2}_*i_2^*2^{m-2}\xi_2\dots \xi_m \\
                            &= \xi_1 2^{m-2}\xi_2\dots \xi_m
                            + \xi_1 2^{m-2}\xi_2\dots \xi_m \\
                            &= 2^{m-1}\xi_1\dots \xi_m
\end{align*}
and Lemma~\ref{lem:identify-class} is proved.
\end{proof}

This proves that $c_{2}$, \dots,~$c_{n}$, $y_{m}$ generate $\A_{\SO}$.

\subsubsection*{Step 2: the relations are satisfied}

The fact that $2c_{i} = 0$ when $i$ is odd follows immediately, as for $\O$, from the fact that $V$ is self-dual.

To prove that $y_mc_{i}=0$, it is sufficient to show that $c_m(L_{m})c_{i}=0$ in $\A(Q_1)$, for any vector bundle $E$ on $X$, with an $\SO$ structure, as $y_mc_{i}={f_{1}}_*(s \cdot \c_m(L_{m})c_{i})$. But on $Q_1$ there is an exact sequence of vector bundles
   \[
   0 \arr L_{m} \arr f^*E \arr L_{m}^{\vee} \arr 0
   \]
so the total Chern class $\c({f_{1}}^*E)$ is $\c(L_{m})\c(L_{m}^{\vee})$ and $\c_{i}(f^*E)=0$ when $i$ is odd.

Finally, the normal bundle $N$ of $C$ in $\AA^n \setminus \{0\}$ is trivial, since the ideal of $C$ is generated by an invariant function on $\An -\0$, so
\begin{align*}
   y_m^2 &= i_*y_{m-1}\cdot i_*y_{m-1}\\
   &= i_*(y_{m-1}\cdot i^*i_*y_{m-1})\\
   &= i_*(y_{m-1}^2\cdot \c_1(N))\\
   &=0
\end{align*}
in $\A_{\SO}(\AA^n \setminus \{0\}) = \A_{\SO} /(c_n)$, by the projection formula and the self-intersection formula. Hence there is an integer $d$ such that $y_m^2 = dc_n$; we will compute $d$ once again by restricting to a maximal torus. By Lemma~\ref{lem:restrict-torus} we have
\begin{align*}
   y_m^2 &= 2^{2m-2} \xi_1^2\dots \xi^2_m\\
   &= 2^{n-2}(-1)^m c_n \in \A[n]_{T_{m}};
\end{align*}
hence, since $c_{n}$ is not a torsion element of $\A_{T_{m}}$, we get that $d = 2^{n-2}$, as claimed.

\subsubsection*{Step 3: the relations suffice}

Consider the ideal $J$ in the polynomial ring $\ZZ[x_{2}, \dots , x_{n},
z]$ generated by the polynomials $z^2 - (-1)^m 2^{n-2} x_{n}$,
$2x_{\mathrm{odd}}$, $z x_{\mathrm{odd}}$. Let $P \in \ZZ[x_{2}, \dots ,
x_{n}, z]$ a homogeneous polynomial such that
   \[
   P(c_{2}, \dots,c_{n}, y_{m}) = 0;
   \]
we need to show that $P$ is in $J$.

By modifying $P$ by an element of $J$, we may assume that it is of the
form $Q_{1} + zQ_{2} + R$, where $Q_{1}$ and $Q_{2}$ are polynomials in
the even $x_{i}$, while $R$ is a polynomial in the $x_{i}$ with
coefficients that are all $0$ or $1$, and all of whose non-zero monomial
contain some $x_{i}$ with $i$ odd.

The odd $c_{i}$ restrict to $0$ in $\A_{T_{m}}$, while $c_{2j}$
restricts to the $j\th$ symmetric function $s_{j}$ of $-\xi_{1}^{2}$,
\dots,~$-\xi_{m}^{2}$; also, $y_{m}$ restricts to $\xi_{1}\dots
\xi_{m}$. Hence $P(c_{2}, \dots, c_{m}, y_{m}) = 0$ restricts to
$Q_{1}(s_{2}, s_{4}, \dots) + \xi_{1}\dots \xi_{m} Q_{2}(s_{2}, s_{4},
\dots)$; and this is easily seen to imply that $Q_{1} = Q_{2} = 0$.

Hence $P$ is a polynomial in $x_{2}$, \dots,~$x_{n}$, all of whose coefficients are $0$ or $1$. Now consider the basis $e'_{1}$, \dots,~$e'_{n}$ of $V$, and the subgroup $\mmu^{n} \subseteq \O$ considered in the previous section, consisting of linear transformations that take each $e'_{i}$ into $e'_{i}$ or $-e'_{i}$. The subgroup $\Gamma_{n} \eqdef \mmu^{n} \cap \SO$ consists of the elements $(\epsilon_{1}, \dots, \epsilon_{n})$ of $\mmu^{n}$ such that $\epsilon_{1}\dots \epsilon_{n} = 1$ in $\mmu$. The group $\Gamma_{n}$ is isomorphic to $\mmu^{n-1}$; if we call $\eta_{i} \in \A[1]_{\Gamma_{n}}$ the first Chern class of the restriction to $\Gamma_{n}$ of the $i\th$ projection $\mmu^{n} \arr \mmu \subseteq \gm$, then we have
   \[
   \A_{\Gamma_{n}} = \ZZ[\eta_{1}, \dots, \eta_{n}]/
   (\eta_{1} + \dots + \eta_{n}).
   \]

We have a natural homomorphism $\A_{\Gamma_{n}}\arr \FF_{2}[\eta_{1}, \dots, \eta_{n}]/(\eta_{1} + \dots + \eta_{n})$, which is an isomorphism in positive degree. If we denote by $r_{1}$, \dots,~$r_{n}$ the elementary symmetric functions of the $h_{i}$, we have that $c_{i}$ restricts to the image of $r_{i}$ in $\FF_{2}[\eta_{1}, \dots, \eta_{n}]/(r_{1})$; hence all we need to show is that the images of $r_{2}$, \dots,~$r_{n}$ are algebraically independent in $\FF_{2}[\eta_{1}, \dots, \eta_{n}]/(r_{1})$. But $r_{1}$, \dots,~$r_{n}$ are algebraically independent in $\FF_{2}[\eta_{1}, \dots, \eta_{n}]$, so $r_{2}$, \dots,~$r_{n}$ are algebraically independent in $\FF_{2}[r_{1}, \dots, r_{n}]/(r_{1})$; and the homomorphism
   \[
   \FF_{2}[r_{1}, \dots, r_{n}]/(r_{1}) \arr
   \FF_{2}[\eta_{1}, \dots, \eta_{n}]/(r_{1})
   \]
is injective, because the extension $\FF_{2}[r_{1}, \dots, r_{n}] \subseteq \FF_{2}[\eta_{1}, \dots, \eta_{n}]$ is faithfully flat. This shows that $P = 0$, and completes the proof of the theorem.

\bibliographystyle{amsalpha}
\bibliography{mrabbrev,VistoliRefs}

\end{document}